\documentclass[12pt]{article}

\usepackage{amsthm}

\title{Some remarks about Descartes' rule of signs
}

\author{ Alain Albouy$^{1,2}$,\qquad
Yanning Fu$^{1}$\\ \\
 $^{1}$ Purple Mountain Observatory\\
Chinese Academy of Sciences\\
2 West Beijing road, Nanjing 210008, China\\
$^{2}$ CNRS-UMR8028, IMCCE\\Observatoire de Paris\\
  77 avenue Denfert-Rochereau, 75014 Paris, France\\
  albouy@imcce.fr, \qquad fyn@pmo.ac.cn}

\date{}

\begin{document}

\maketitle

This preprint is now published. DOI: 10.4171/EM/262

Reference: {\it Elemente der Mathematik}, 69 (2014), pp.\ 186--194

\bigskip

\noindent{\bf Introduction}

\noindent What can we deduce about the roots of a real polynomial in one variable by only considering the signs of its coefficients?

\noindent If we restrict our attention to the positive roots, the answer is: nothing more than what the classical Descartes rule of signs states. This is a precise mathematical claim that needs a proof. We give it in
Section 3, after some due preliminaries.

\noindent If we are interested in the negative roots also, and concentrate on the combined numbers of positive and negative roots, there is nothing as simple to tell. In Section 5, we show how the situation complicates with the degree, by giving a complete study up to degree 6.

\noindent Section 4 is a little interlude about discriminants and trinomials.

\bigskip
\centerline{\bf 1. Coefficients and positive roots.}
\bigskip

It is convenient to consider that an analytic function in one variable has a {\it
multiset} of roots. A multiset is defined just as a set, except that  identical elements
are allowed. The {\it cardinality} of a multiset is the number of its elements. The {\it multiplicity} of an element is the cardinality of the sub-multiset of the elements identical to it.

Consider expressions of the type
\begin{equation}\label{Y}
Y=a_0x^{\alpha_0}+a_1x^{\alpha_1}+\cdots+a_nx^{\alpha_n},
\end{equation}
with real exponents $\alpha_0<\alpha_1<\cdots<\alpha_n$ and non-zero real coefficients $a_0,\dots, a_n$.
The {\it sequence of signs} of an expression $(\ref{Y})$ is the ordered list $\sigma_0,\dots,\sigma_n$ of the signs $\sigma_i=|a_i|^{-1}a_i$. A {\it variation} occurs in such a sequence when $\sigma_i=-\sigma_{i+1}$.  For example the sequence $1,-1,1,1$ has two variations.

{\bf Proposition 1.} {\it  On one hand, let us prescribe the number of terms, the exponents and the signs of the coefficients in an expression $(\ref{Y})$, and denote by $k$ the number of variations of the sign sequence. On the other hand, let us give arbitrarily a multiset of positive numbers, and call $p$ its cardinality. There exists an expression $(\ref{Y})$, as prescribed, whose multiset of positive roots is the given one, if and only if $k-p$ is even and non-negative.}

This statement gives an ``if and only if" form to the classical Descartes rule of signs. It continues previous works. In 1998, Anderson, Jackson and Sitharam \cite{AJS} proposed examples of polynomials with any sequence of signs, and with any number $p$ of positive roots, provided that $k-p$ is even and non-negative.  In 1999, Grabiner \cite{Gra} extended this result by giving also examples of polynomials with missing terms (i.e.\ zero coefficients for some of the intermediate integer powers of $x$). Our Proposition 1 strengthens this statement, by extending these previous results to real exponents, and by showing that when a given sequence of signs allows $p$ positive roots, it indeed allows any multiset of positive roots of cardinality $p$. Surprisingly the proof, as presented below, remains very elementary.

\bigskip
\centerline{\bf 2. Proof of the ``only if" part.}
\bigskip

This proof is well-known and we give it only for completeness. We assume that the given multiset is the multiset of positive roots of $Y$.

(i) The parity of $k$ decides if the sign of $a_0$ is or is not the sign of $a_n$. The parity of $p$ decides if the sign of $Y$ is or is not the same at the two ends, i.e.\ at $0^+$ and at $+\infty$. Obviously these boundary signs are respectively the signs of $a_0$ and of $a_n$. Thus $k-p$ is even.

(ii) Let $m$ and $m'$ be the cardinalities of the respective multisets of positive roots of a function $f$ and of its derivative $f'$. Then $m\leq m'+1$.

(iii) One multiplies $Y$ by $x^{-\alpha}$, computes the derivative in $x$ and observes the resulting sequence of signs. The signs corresponding to the exponents $\alpha_i<\alpha$ are changed, while those corresponding to $\alpha_i>\alpha$ are not changed. We can delete any given variation of sign without touching the others by choosing
$\alpha$ in some interval of the form $]\alpha_i,\alpha_{i+1}[$.

To conclude, we prove that $p\leq k$ by an induction on $k$. If $k=0$ then $p=0$.  If $Y$ has $k$  variations, argument (iii) shows that $(x^{-\alpha}Y)'$ has $k-1$ variations for some choice of $\alpha$. Argument (ii) gives the required estimate on the respective numbers of roots.\qed

These arguments are actually older than usually thought. Argument (ii) was given in 1741 by de Gua  and again in 1798 by Lagrange. Argument (iii) and the statement of the ``only if" part are given by de Gua in the case of integer exponents, by Laguerre in the case of real exponents\footnote{See \cite{Gua} pp.\ 89--92, \cite{Lag} p.\ 195,  \cite{Lagu}. Jensen \cite{Jen} noticed the similarity of de Gua's and Laguerre's arguments. Instead of derivation, de Gua used Johan Hudde's operation, which is described as a term by term multiplication of the polynomial and an arithmetic progression. The modern reader will simply compare this with the operation $Y\mapsto x^{-m+1}(x^mY)'$.}. Laguerre is explicit about the parity conclusion (i), which is obvious. If the statement of the ``only if" part omits this conclusion, then it remains correct if the roots are counted without multiplicity\footnote{Authors of the 18th century used to count the roots with their multiplicity and to omit the parity conclusion. Incidentally, many authors following Cajori \cite{Caj} attribute the parity conclusion to Gauss \cite{Gau} in 1828. This is very strange as firstly, this conclusion is obvious, secondly, Fourier published it in 1820 (see  \cite{Fou} p.\ 294) and thirdly, it is absent from Gauss' paper.}.

\bigskip
\centerline{\bf 3. Proof of the ``if" part.}
\bigskip

We call $(\sigma_0,\dots,\sigma_n)$ the prescribed sequence of signs (with $k$ variations). We assume that the cardinality $p$ of the given multiset of positive numbers is such that $k-p$ is  non-negative and even.
We look for an expression $(\ref{Y})$ with $\sigma_ia_i>0$ for all
$i$, such that the given multiset is the multiset of its positive roots.
We will be able to find such an expression with the additional constraint that $a_i=a_{i+1}$ when $\sigma_i=\sigma_{i+1}$.
We  collect the terms in $(\ref{Y})$ according to this constraint. We set $i_0=0$ and call
$i_1, i_2,\dots, i_k$ the integers $i_j$ such that $i_{j-1}<i_j$, $\sigma_{i_j-1}\neq\sigma_{i_j}$, $\sigma_{i_k}= \sigma_n$. We call $b_0,b_1,\dots,b_k$ the common values of the coefficients:
$$b_j=a_{i_j}=a_{i_j+1}=\cdots=a_{i_{j+1}-1}.$$
We set $$\varphi_j(x)=x^{\alpha_{i_j}}+x^{\alpha_{i_j+1}}+\cdots+x^{\alpha_{i_{j+1}-1}}.$$
Now we look for a $(b_0,b_1,\dots,b_k)$ such that $\sigma_0 b_0>0$, $b_ib_{i+1}<0$ and the expression
$$Y=b_0\varphi_0(x)+b_1\varphi_1(x)+\cdots+b_k\varphi_k(x)$$
has the required multiset of positive roots.

We begin with the case where all the positive roots $x_1,\dots,x_p$ are simple. The column vector $(b_0,\dots,b_k)$ should be in the kernel of the matrix
$$\Phi=\pmatrix{\varphi_0(x_1)&\varphi_1(x_1)&\cdots&\varphi_k(x_1)\cr
\varphi_0(x_2)&\varphi_1(x_2)&\cdots&\varphi_k(x_2)\cr
&&\cdots&\cr
\varphi_0(x_p)&\varphi_1(x_p)&\cdots&\varphi_k(x_p)\cr}.$$
This $p\times (k+1)$ matrix has a non-trivial $(b_0,\dots,b_k)$ in its kernel since $p\leq k$.

If $p=k$, the ``only if" part shows that the coefficients $b_i$ are non-zero, that their signs alternate and that there is no other positive root. This makes an $(a_0,\dots , a_n)$ such that $Y$ has the prescribed set of positive roots.

If $1\leq p<k$, observe that the $p\times p$ submatrix formed by the first $p$ columns of $\Phi$ is invertible. If its determinant were zero we could find a non-trivial column vector $(c_0,\dots,c_{p-1})$ in its kernel. But $\sum_0^{p-1} c_i\varphi_i(x)$ would have a root at each $x_j$. It would have $p$ roots and at most $p-1$ variations, which contradicts what we proved in the ``only if" part.

So for any given $(b_p,\dots, b_k)$ there is a unique  $(b_0,\dots,b_{p-1})$, computed by Cramer's rule, such that $(b_0,\dots,b_k)$ is in the kernel of $\Phi$. We introduce a real parameter $\epsilon$ and take $b_i(\epsilon)= (-1)^{k-i}\sigma_n\epsilon$, $p\leq i<k$ and $b_k(\epsilon)=\sigma_n$. The unique column vector $\bigl(b_0(\epsilon),\dots, b_{p-1}(\epsilon),b_p(\epsilon),\dots,b_k(\epsilon)\bigr)$ in the kernel of $\Phi$ is such that each of its $k+1$ components is a polynomial in $\epsilon$, of first degree.

Set $Y_\epsilon(x)=b_0(\epsilon)\varphi_0(x)+\cdots+b_k(\epsilon)\varphi_k(x)$. We have $Y_\epsilon(x_j)=0$ for $j=1,\dots,p$.
Since $Y_0(x)=b_0(0)\varphi_0(x)+\cdots+b_{p-1}(0)\varphi_{p-1}(x)+b_k(0)\varphi_k(x)$, the $p+1$ coefficients in this expression are non-zero and have alternating signs, according to the ``only if" part, which at the same time establishes that $Y_0$ cannot have any other positive root than the $x_j$'s. As $k-p$ is even and $b_k$ has the prescribed sign, the $b_i(0)$'s, $i<p$, have the prescribed signs, and the same is true of $b_i(\epsilon)$ for any small enough $\epsilon$. Now, the $b_i(\epsilon)$'s, $p\leq i<k$, have the prescribed signs for any positive $\epsilon$. Thus $Y_\epsilon$ has the prescribed sign sequence for any sufficiently small positive $\epsilon$.

Consider $Y_\epsilon/Y_0$. This quotient is analytic in $x\in ]0,+\infty[$, as shown by expanding numerator and denominator in Taylor series around  $x=x_j$ and by simplifying the factors $(x-x_j)$.
As $\epsilon\to 0$, $Y_\epsilon/Y_0\to 1$ for any $x\in]0,+\infty[$. This convergence is indeed uniform, the behavior of the fraction at zero and at infinity being controlled by the leading terms of the numerator and the denominator. Thus $Y_\epsilon$ has no other positive root for a sufficiently small $\epsilon$. The problem is solved in the case of simple roots.

In the cases with multiple roots the construction and the proofs are exactly the same, except that $\Phi$ has e.g.\ the row $\bigl(\varphi_0'(x_1),\dots,\varphi_k'(x_1)\bigr)$ if $x_1$ is a multiple root, the row $\bigl(\varphi_0''(x_1),\dots,\varphi_k''(x_1)\bigr)$ if it is at least a triple root, etc.

If finally $p=0$, $k$ is even and we can build a $Y$ without positive root. We simply choose a sufficiently small positive $\epsilon$ and take $b_0(\epsilon)=\sigma_0$, $b_k(\epsilon)=\sigma_0$ and $b_i(\epsilon)=(-1)^i\epsilon\sigma_0$.\qed

\bigskip
\centerline{\bf 4. Further information about trinomials.}
\nobreak
\bigskip
Etymologically a trinomial is simply the sum of three monomials, and there is no restriction on the degree. Trinomials with unprescribed degree were studied early in the history of algebraic equations (see \cite{Ste1}, pp.\ 11 and 24). They are natural objects in the context of Laguerre's extension of Descartes' rule to real exponents, as well as in the context of Khovansky's theory of fewnomials \cite{Kho}, and are consequently the object of recent studies. In 2002, Haas \cite{Haa}, Li, Rojas \& Wang \cite{LRW} proved that the optimal upper bound on the number of roots in the positive quadrant of a system of two trinomials in two variables is five. 

We were not able to find the following elegant and elementary formula in any of these old or recent studies.

{\bf Proposition 2.} {\it A trinomial $a x^\alpha+bx^\beta+c x^\gamma$, $\alpha<\beta<\gamma$, $a>0$, $c>0$, $b<0$,
is positive on $]0,\infty[$ if and only if
$$\Bigl(\frac{a}{\gamma-\beta}\Bigr)^{\gamma-\beta}\Bigl(\frac{b}{\alpha-\gamma}\Bigr)^{\alpha-\gamma}\Bigl(\frac{c}{\beta-\alpha}\Bigr)^{\beta-\alpha}>1.$$}

{\bf Remark.} The discriminant of $a x^\alpha+bx^\beta+c x^\gamma$ is $b^2-4ac$ if $(\alpha,\beta,\gamma)=(0,1,2)$. It is $-c(4b^3+27ca^2)$ if $(\alpha,\beta,\gamma)=(0,1,3)$.
It is $c^2(-27b^4+256ca^3)$ if $(\alpha,\beta,\gamma)=(0,1,4)$ and $16ac(4ca-b^2)^2$ if
$(\alpha,\beta,\gamma)=(0,2,4)$. These formulas are quite familiar. The general formula in the proposition gives in each case the main factor.

{\bf Proof.} The formula may be obtained by a direct computation, where there occur unexpected simplifications. We will present a short matrix algebra argument, basically the same as in the previous section. Consider the trinomial $Y=ax^{\alpha}+bx^{\beta}+cx^{\gamma}$. We determine $(a,b,c)$ such that $Y$ has a double root at $x_1$. We find that $(a,b,c)$ should be in the kernel of
$$\Phi=\pmatrix{x_1^\alpha&x_1^\beta&x_1^\gamma\cr
\alpha x_1^{\alpha-1}&\beta x_1^{\beta-1}&\gamma x_1^{\gamma-1}},$$
which means
$$\frac{ax_1^\alpha}{\gamma-\beta}=\frac{bx_1^\beta}{\alpha-\gamma}=\frac{cx_1^\gamma}{\beta-\alpha}.$$
We set $$A=\frac{a}{\gamma-\beta},\quad B=\frac{b}{\alpha-\gamma},\quad C=\frac{c}{\beta-\alpha},$$
which are positive numbers according to the hypotheses of the proposition, and continue the computation by eliminating $x_1$, which gives the expected condition $A^{\gamma-\beta}B^{\alpha-\gamma}C^{\beta-\alpha}=1$ for a double root. Now a mere study of $A^{\gamma-\beta}B^{\alpha-\gamma}C^{\beta-\alpha}$ as a function of $b$ gives the proposition.\qed

For completeness, note that under the condition $A^{\gamma-\beta}B^{\alpha-\gamma}C^{\beta-\alpha}=1$ the double root is located at
$$x_1=\Bigl(\frac{C}{B}\Bigr)^{1/(\beta-\gamma)}=\Bigl(\frac{A}{C}\Bigr)^{1/(\gamma-\alpha)}=\Bigl(\frac{B}{A}\Bigr)^{1/(\alpha-\beta)}.$$

\bigskip
\centerline{\bf  5. Combining positive and negative roots.}
\nobreak
\bigskip

If we just focus on the positive roots of a given polynomial $Y$, Proposition~1 tells us that there is nothing we can add to Descartes' rule of signs. All that can be deduced from
the sign sequence of $Y$ is the upper bound and the parity of $P$,  the cardinality of the multiset of positive roots. By changing $x$ to $-x$, the same can be said on $N$, the cardinality of the multiset of negative roots of $Y$.  But let us consider the constraints on  $(P,N)$.

Consider a polynomial $Y(x)$ of degree 4 with sign sequence
$+,-,-,-,+$. By Descartes' rule  $P=0$ or 2. If we change $x$ to $-x$ the sign sequence becomes $+,+,-,+,+$. Thus $N=0$ or 2. Grabiner
\cite{Gra} points out that $(P,N)=(0,2)$ is {\it impossible} for such a $Y(x)$.

The proof does not require any computation. Because the constant term $Y(0)$ is
positive, $P=0$ implies $Y>0$ for $x>0$. But the odd part of $Y$ is negative
when $x>0$. So the even part is positive. For $x<0$ the odd and the even parts are then
positive. We must have $N=0$.\qed

Consider an expression $(\ref{Y})$ as in Proposition 1, where the $\alpha_i$'s are
non-negative integers, i.e.\ $Y(x)$ is a polynomial. Given the sequence of signs
$\sigma_0,\dots,\sigma_n$, Grabiner proves\footnote{The idea is that small terms can be neglected before applying Descartes' rule. Estimates on how small these terms should be are discussed in \cite{Mar}. It is interesting to compare these estimates to the first Lemma in \cite{Gua}.} that all the $(P,N)$'s produced as follows are {\it possible}.
One chooses any subset of $\{1,\dots, n-1\}$ and ``erases" the corresponding
$\sigma_i$'s from the sequence of signs. The number of variations in the
resulting sequence gives $P$. One then considers the modified sequence of signs
$\sigma_0(-1)^{\alpha_0},\sigma_1(-1)^{\alpha_1},\dots,\sigma_n(-1)^{\alpha_n}$
and erases from it the signs with index in the same subset. The new number of
variations gives $N$.

Let us call this construction of possible $(P,N)$'s Grabiner's erasing term rule.
Grabiner observes that it gives all the possible $(P,N)$'s for polynomials
up to degree 4, and conjectures that the same is true for higher degrees.

Here is a counterexample. Consider a polynomial of 5th degree with sign
sequence  $+,+,-,+,+,-$. Here $(P,N)=(3,0)$ is possible, as shown by
$$ \scriptstyle{(20 + 37 x + 18 x^2) (1-x) (2-x)
(3-x) =  120 +2 x -179 x^2+ 4 x^3+71 x^4-18 x^5.}$$ This
combination is not obtained by Grabiner's erasing term rule. Indeed, after changing
$x$ in $-x$ the sequence becomes $+,-,-,-,+,+$. As $N=0$ we should erase terms
and obtain a sequence without variation. We should erase either the three $-$,
or the three $-$ and the internal $+$. But erasing the corresponding signs in
$+,+,-,+,+,-$ gives in both cases only one variation, while $P=3$.

In contrast $(P,N)=(3,0)$ is impossible for the sequence $+,+,-,+,-,-$.
Descartes' rule predicts all the impossibilities for degree 5 polynomials without gaps, except this one and its trivial analogues.

To prove this impossibility, let us specify the ordering of sign
sequences by associating the first sign to the constant term.
We notice that the odd part ${\cal O}$ has sequence of signs $+,+,-$. Its derivative of
course has the same sequence and consequently a unique positive root $r$.
On $]0,r[$, ${\cal O}$ is
positive, and for $x>r$, ${\cal O}$ is decreasing.
The even part ${\cal E}$ has sequence of signs $+,-,-$, thus it decreases for
$x>0$. Our polynomial ${\cal E}+{\cal O}$ decreases and has at most 1 root on $]r,+\infty[$.
As $N=0$, ${\cal E}-{\cal O}$ is
positive on $]0,+\infty[$.  On $]0,r[$, ${\cal O}$ and thus ${\cal E}+{\cal O}$ are positive. There is one positive root on $]r,+\infty[$, thus $P=3$ is impossible.\qed

Let us pass to degree 6 polynomials. The exhaustive list of non-Descartes impossibilities is, up to trivial transformations:
$${+,+,-,+,-,+,+} \hbox{ is incompatible with }(P,N)=(2,0)\hbox{ or }
(4,0)$$
$${+,+,+,+,-,+,+} \hbox{ is incompatible with }(P,N)=(2,0)$$
$${+,+,-,-,-,-,+} \hbox{ is incompatible with }(P,N)=(0,4)$$
The argument we gave in Grabiner's example proves the first and second
statements.

To prove the third
statement, we write a polynomial with $(P,N)=(0,4)$ as
$$
p = (c - b x + x^2) (x + x_1) (x + x_2) (x + x_3) (x + x_4),
$$
where $x_i>0$, $i=1,\dots,4$,
\begin{equation}\label{0}
4c>b^2>0.
\end{equation}
Expanding,
$
p = a_0+a_1x+\cdots+a_5x^5+x^6,
$
with in particular
\begin{equation}\label{a2}
a_2= c \beta - b \gamma + \delta<0, \\
\end{equation}
\begin{equation}\label{a5}
a_5=-b + \alpha<0,
\end{equation}
the positive numbers $\alpha, \beta,\gamma,\delta$ being defined
as
\begin{equation}\label{xi}
\begin{array}{l}
\alpha = x_1 + x_2 + x_3 + x_4 \\
\beta= x_1 x_2 + x_1 x_3 + x_2 x_3 + x_1 x_4 + x_2 x_4 + x_3 x_4\\
\gamma = x_1 x_2 x_3 + x_1 x_2 x_4 + x_1 x_3 x_4 + x_2 x_3 x_4 \\
\delta = x_1 x_2 x_3 x_4 \end{array}
\end{equation}
By $(\ref{0})$ and $(\ref{a5})$,
\begin{equation}\label{I}
\sqrt{c}>b/2>\alpha/2>0,
\end{equation}
and, by $(\ref{0})$ and $(\ref{a2})$,
$$
2\sqrt{c}>b>\frac{c \beta+ \delta}{\gamma},
$$
which gives
\begin{equation}\label{tmp2}
\beta \sqrt{c}^2-2\gamma \sqrt{c}+\delta<0.
\end{equation}
We expand the discriminant $4\gamma^2-4 \beta \delta$ of this expression using $(\ref{xi})$. All the terms are positive. Now $(\ref{tmp2})$ implies in particular that
$$
\sqrt{c}<\frac{\gamma+\sqrt{\gamma^2-\beta\delta}}{\beta}.
$$
Combining with $(\ref{I})$ gives
$
\beta\alpha/2<\beta\sqrt{c}<\gamma+\sqrt{\gamma^2-\beta\delta}
$
or
\begin{equation}\label{A}
\alpha\beta/2-\gamma<\sqrt{\gamma^2-\beta\delta}.
\end{equation}
After expanding and
canceling, we see that all the terms of $\alpha\beta/2-\gamma$ are positive. We
square both sides of $(\ref{A})$. This
gives us
$
\alpha^2 \beta-4\alpha\gamma+4\delta<0.
$
But the left hand side, after expanding and
canceling, has only positive terms. This is a contradiction.\qed

To check that our list of impossibilities is complete at degree 6, it is enough to consider the possibilities established by Grabiner's erasing
term rule and the two following polynomials, which after changing $x$ into
$-x$, $x$ into $1/x$ or $Y$ into $-Y$ provide  examples for all the other possibilities:
$$\scriptstyle{(31+11x+x^2)(1-x)(2-x)(3-x)(4-x)=744-1286x+559x^2+25x^3-44x^4+x^5+x^6,}$$
$$\scriptstyle{(9+8x+2x^2)(1-x)(2-x)(3-x)(4-x)=216-258x-37x^2+90x^3-x^4-12x^5+2x^6.}$$
These possibilities and impossibilities do not seem to organize themselves in simple classes.

{\bf Acknowledgements.} We wish to thank Alain Chenciner, David Grabiner, Dieter Schmidt and the referee for their precious remarks.  This work is supported by Chinese Academy of Sciences visiting professorship for senior international scientists grant no.\ 2009J2-11, and NSFC grant no.\ 10833001.

\end{document}